\documentclass[12pt]{article}
\usepackage[a4paper, left=3.3cm, right=3.3cm, top=3cm]{geometry}
\usepackage[utf8]{inputenc}
\usepackage{amsmath}
\usepackage{amsthm}
\usepackage{amsfonts}
\usepackage{amssymb}
\usepackage{graphicx}
\usepackage{color}

\usepackage{epstopdf}
\usepackage[caption=false]{subfig}
\usepackage{color}

\usepackage[numbers,sort&compress]{natbib}
\bibpunct[, ]{[}{]}{,}{n}{,}{,}
\renewcommand\bibfont{\fontsize{10}{12}\selectfont}
\makeatletter
\def\NAT@def@citea{\def\@citea{\NAT@separator}}
\makeatother

\theoremstyle{plain}
\newtheorem{theorem}{Theorem}[section]
\newtheorem{lemma}[theorem]{Lemma}

\newtheorem{proposition}[theorem]{Proposition}

\theoremstyle{definition}
\newtheorem{definition}[theorem]{Definition}
\newtheorem{example}[theorem]{Example}

\theoremstyle{remark}
\newtheorem{remark}{Remark}

\newcommand{\R}{\mathbb{R}}
\newcommand{\Z}{\mathbb{Z}}
\newcommand{\N}{\mathbb{N}}
\newcommand{\Oe}{\Omega_\epsilon}
\newcommand{\OseM}{\Omega_{*,\epsilon}^M}
\newcommand{\OeM}{\Omega_{\epsilon}^M}
\newcommand{\Ssepm}{S_{*,\epsilon}^\pm}

\begin{document}
	
	
	\title{Effective transmission conditions for reaction-diffusion processes in domains separated by thin channels}
	
	\author{Apratim Bhattacharya\thanks{Department Mathematik, FAU Erlangen-N\"urnberg, Cauerstr. 11, 91058 Erlangen, Germany}, Markus Gahn\thanks{IWR, Universit\"at Heidelberg, INF 205, 69112 Heidelberg, Germany}, Maria Neuss-Radu\thanks{Department Mathematik, FAU Erlangen-N\"urnberg, Cauerstr. 11, 91058 Erlangen, Germany}
	}
	
	\maketitle
	
	\begin{abstract}
		We consider a reaction--diffusion equation in a domain consisting of two bulk regions connected via small channels periodically distributed within a thin layer. The height and the thickness of the channels are of order $\epsilon$, and the equation inside the layer depends on the parameter $\epsilon$ and an additional parameter $\gamma \in [-1,1)$, which describes the size of the diffusion in the layer. We derive effective models for the limit $\epsilon \to 0 $, when the channel-domain is replaced by an interface between the two bulk-domains. 
	\end{abstract}
	
	\noindent\textbf{Keywords:}
	Array of channels; homogenization; two-scale convergence; reaction-diffusion equation; effective transmission conditions
	
	\noindent\textbf{AMS Subject Classification:} 35K57; 35B27, 35Q92
	
	\section{Introduction}
	
	In this paper, we consider a reaction--diffusion equation in a domain $\Omega$ consisting of two bulk regions $\Omega_\epsilon^+$ and $\Omega_\epsilon^-$ which are connected via small periodically distributed channels. The height and the thickness of the channels are of order $\epsilon$, where the parameter $\epsilon$ is small compared to the size of the bulk-domains. The equation in the channels' domain depends on $\epsilon$, the diffusion coefficients having the size $\epsilon^\gamma$ with $\gamma \in [-1,1)$. On the lateral boundary of the channels, we consider a Neumann-boundary condition describing e.g., surface reactions. The aim of this paper is the derivation of an effective model in the limit $\epsilon \to 0$. This requires asymptotic techniques which combine the classical approach of homogenization for microscopic structures with a singular limit approach, in the context of a thin layer perforated by disconnected channels. In the limit $\epsilon \to 0$, the channel-domain is replaced by an interface $\Sigma$ between the bulk-domains $\Omega^+$ and $\Omega^-$, and the main difficulty lies in the derivation of effective interface laws on $\Sigma$. Here, we extend reaction-diffusion problems and corresponding techniques considered previously in \cite{Neu07, Gah17, Gah18} to the situation of channels. We consider the case $\gamma \in [-1,1)$, i.e., moderate and high diffusion inside the channels' domain. In the limit $\epsilon \to 0$, at the interface $\Sigma$, we obtain the continuity of the concentrations and an ordinary differential equation for the limit concentration in the layer, involving the jump of the normal fluxes of the solutions in $\Omega^\pm$ across $\Sigma$. The critical case $\gamma = 1$ is treated in \cite{Gah20} and leads to different effective transmission conditions at $\Sigma$, namely a jump in the homogenized solution and its normal fluxes depending on solutions to local problems in the standard channel in every macroscopic variable $\bar x \in \Sigma$. We consider solutions with low regularity with respect to the time variable in order to keep the assumptions on the data quite general. 
	
	First homogenization results for problems with a geometrical framework related to our setting were given in \cite{San82}. Contributions to the homogenization of the Laplace equation in domains connected by thin channels have been given in \cite{Del87, Sha01, Yab04, Ono06}, where the asymptotic behavior of the solution is investigated for different ratio of the thickness of the layer and the radius of the cylindrical channels. The homogenization of an elliptic Steklov type spectral problem in domains connected by thin channels was considered in \cite{Ami16, Gad18}. Reaction-diffusion problems through channels with low conducting properties were considered in \cite{Gah12}, and is also topic of ongoing work, see \cite{Gah20}. The more challenging problem concerning the ion transport through channels of biological membranes was announced in \cite{Neu05}. 
	
	This paper is organized as follows: In Section \ref{Setting_of_Problem}, we define the microscopic geometry, introduce the microscopic model, and give its variational formulation with respect to function spaces adapted to the scaling of the model. In Section \ref{A_priori_estimates}, we derive estimates for the microscopic solutions necessary for the derivation of compactness results, especially with respect to the two-scale convergence for channels. The definition of the latter together with the compactness results are given in Section \ref{Two_scale_convergence}. In Section \ref{Convergence_results} the convergence results for the microscopic solutions are derived. Especially the case $\gamma = -1$, when the gradient of the solution in the channels converges to zero in the two-scale sense, requires refined techniques adapted to the microscopic geometry. In Section \ref{Derivation_macroscopic_model} the homogenized model, including the effective interface conditions is derived.
	
	\section{Setting of the Problem}
	\label{Setting_of_Problem}
	
	Let $\epsilon > 0$ be a sequence of positive numbers tending to zero such that $\frac{1}{\epsilon} \in \N$ and let $H>0$ be a fixed real number. Let $n\in \N, n\geq2$. For $x\in \R^n$, we write $x=(\bar x, x_n) \in \R^{n-1}\times \R$.
	Let $\Omega$ be a subset of $\R^n$ defined as
	$$
	\Omega = \Sigma \times (-H,H), 
	$$
	where $\Sigma \subset \R^{n-1}$ is a rectangle in $\R^{n-1}$, i.e.~$\Sigma = \prod_{i=1}^{n-1} (a_i,b_i)$ with $a_i, b_i \in \Z, a_i <b_i.$ 
	The meaningful case for applications is $n=3$. However, many of our statements can be formulated in $\mathbb{R}^n$ for $n\geq 2$. Therefore, we introduce our notations for a general dimension $n$.
	
\begin{figure}
		\centering
		\subfloat[Microscopic domain $\Omega_\epsilon$ for the case $\epsilon= \frac{1}{3}$ and $n=3$.]{%
			\resizebox*{5cm}{!}{\includegraphics{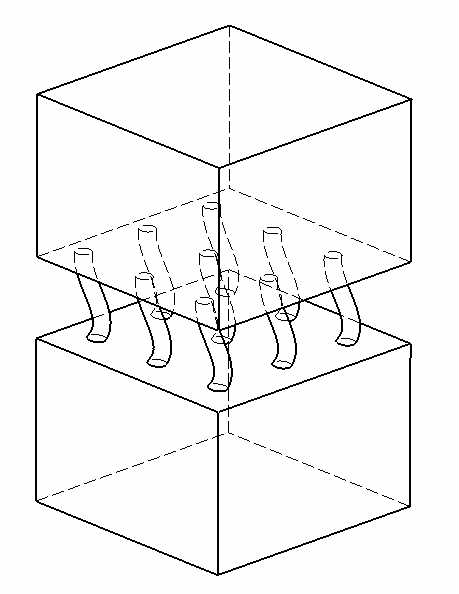}}}\hspace{5pt}
		\subfloat[Standard channel-domain $Z^*$ in the standard cell $Z$.]{%
			\resizebox*{5cm}{!}{\includegraphics{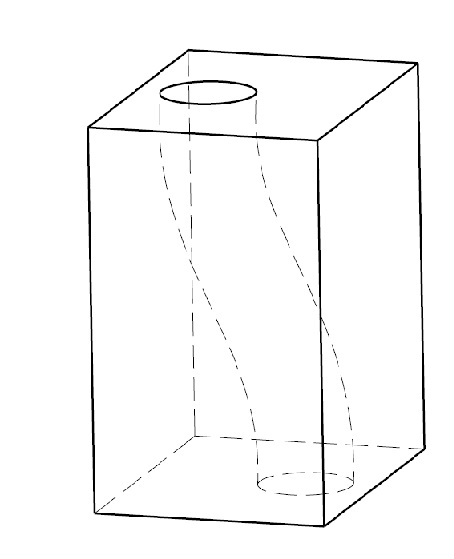}}}
		\caption{Example of a microscopic domain generated by periodically distributed channels.} \label{Microscopic_Domain}
	\end{figure}

	We consider the domain $\Omega_\epsilon \subset \Omega$ consisting of three subdomains: the bulk regions $\Omega_\epsilon^+$ and $\Omega_\epsilon^-$ which are connected by channels periodically distributed within a thin layer constituting the domain $\OseM$, see Figure \ref{Microscopic_Domain} (a). 
	The bulk regions are given by: 
	\begin{equation*}
	\Omega_\epsilon^+ =  \Sigma \times (\epsilon,H), \quad \Omega_\epsilon^- = \Sigma \times (-H,-\epsilon).
	\end{equation*}
	Furthermore, we denote by 
	\begin{equation*}
	S_\epsilon^+ = \Sigma \times \{\epsilon\}, \quad S_\epsilon^- =  \Sigma \times \{-\epsilon\},.
	\end{equation*}
	The thin layer separating the two bulk-domains is given by
	$$
	\OeM = \Sigma \times (-\epsilon, \epsilon).
	$$
	To define the channels, which are periodically distributed within $\OeM$, we first define the standard cell
	$$
	Z:=Y \times (-1,1):= (0,1)^{n-1} \times (-1,1)
	$$ 
	with the upper and lower boundaries
	$$
	S^\pm := \{y= (\bar{y},y_n)\in \R^n: \overline{y} \in Y, y_n =\pm 1\}.
	$$
	Let $Z^* \subset Z$ be a connected and open Lipschitz domain representing the standard channel-domain, such that
	\begin{equation*}
	S_*^\pm :=\{y \in \partial Z^* : y_n=\pm 1\} \neq \emptyset
	\end{equation*}
	is a Lipschitz domain in $\R^{n-1}$ with positive measure, see Figure \ref{Microscopic_Domain} (b). Let the lateral boundary of the standard channel be denoted by 
	$$
	N:= \partial  Z^* \setminus \big(\overline{S_*^+} \cup \overline{S_*^-}\big).
	$$
	We assume that $N$ has a finite distance to the lateral boundary  $\partial Z \setminus (S^+ \cup S^-)$ of the standard cell $Z$. 
	The domain consisting of the channels is then given by
	\begin{equation*} 
	\OseM := \bigcup\limits_{\bar k \in  I_\epsilon}  \epsilon Z_{\bar k}^* : = \bigcup\limits_{\bar k \in  I_\epsilon} \epsilon \big(Z^*+(\bar k,0)\big),
	\end{equation*}  
	where $I_\epsilon= \{ \bar k\in \Z^{n-1}: \epsilon \big(Z+(\bar k,0)\big) \subset \Omega_\epsilon^M \}.$ We have $|I_\epsilon|= \frac{|\Sigma|}{\epsilon ^{n-1}}$. The interfaces between the channel-domain  and the bulk-domains are defined by 
	\begin{equation*}
	\Ssepm := \bigcup\limits_{\bar{k} \in  I_\epsilon} \epsilon \big(S_*^\pm + (\bar k,0)\big).
	\end{equation*}  
	The domain $\Oe$ is thus defined by
	$$
	\Oe= \Omega_\epsilon^+ \cup \Omega_\epsilon^-\cup \OseM \cup S_{*, \epsilon}^+ \cup S_{*, \epsilon}^-.
	$$ 
	We assume $\Omega_\epsilon$ to be Lipschitz. The boundary $\partial \Omega_\epsilon$ can be decomposed into two disjoint subsets $\partial \Omega_ \epsilon=\overline{N_\epsilon}\cup \partial_{N} \Omega_\epsilon$, where $N_\epsilon$ represents the lateral boundary of the channels 
	\begin{equation*}
	N_\epsilon  := \bigcup\limits_{\bar k \in  I_\epsilon} \epsilon \big(N + (\bar k,0)\big).
	\end{equation*}	
	Let $\nu$ denote the outward unit normal at the boundary of $\Omega_\epsilon$. At the interfaces $S_\epsilon^\pm$, we denote by $\nu$ the outward unit normal with respect to $\Omega_\epsilon^\pm$. For a function defined on $\Omega_\epsilon$, we add superscripts $+,-,M$ to denote its restriction to the sub-domains $\Omega_\epsilon^+, \Omega_\epsilon^-$ and $\OseM$ respectively. Finally, we define the domains
	\begin{equation*}
	\Omega^+ := \Sigma \times (0, H), \quad \Omega^-  := \Sigma \times (-H, 0). 
	\end{equation*}
	which are separated by the interface $\Sigma$. We emphasize that, by abuse of notation, we will identify $\Sigma \subset \R^{n-1}$ with $\Sigma \times\{0\} \subset \R^{n}.$ The outer normal to $\Omega^\pm$ is denoted by $\nu^\pm$.
	
	\subsection{The microscopic model}
	
	Let $\gamma \in [-1,1)$. In the domain $\Omega_\epsilon$, we study the following reaction-diffusion equations for the unknown function $u_\epsilon = (u_\epsilon^+, u_\epsilon^M, u_\epsilon^-):(0,T)\times \Omega_\epsilon \rightarrow \R$:
	\begin{subequations}\label{microscopic_model}
		\begin{equation}\label{micro_reac_diff_eq}
		\begin{aligned} 
		&&\partial_{t} u_\epsilon^\pm -D^\pm \Delta u_\epsilon ^\pm &=f^\pm(t,x) &&\text{in  $(0,T)\times  \Omega_\epsilon ^\pm$},  \\
		&&\frac{1}{\epsilon}\partial_{t}u_\epsilon^M - \nabla\cdot\left(\epsilon^\gamma D^M\left(\frac{x}{\epsilon}\right) \nabla u_\epsilon^M\right) &=  \frac{1}{\epsilon}g_\epsilon(t,x) &&\text{in $(0,T)\times \OseM$},
		\end{aligned}
		\end{equation}
		with the boundary conditions: 
		\begin{equation}\label{micro_boundary_cond}
		\begin{aligned} 
		&&-D^\pm \nabla u_\epsilon^\pm  \cdot\nu &=0&&\text{on $ (0,T) \times \partial_{N}\Omega_\epsilon $},\\
		&&-\epsilon^\gamma D^M\left(\frac{x}{\epsilon}\right)\nabla u_\epsilon^M \cdot \nu&= h_\epsilon(t, x) &&\text{on $(0,T) \times N_\epsilon$},
		\end{aligned}
		\end{equation}
		and the initial conditions:
		\begin{equation} \label{micro_initial_cond_pb_1}
		u_\epsilon(0) =  u_{i, \epsilon} \quad \mbox{in } \Oe.
		\end{equation}
		On the interfaces $\Ssepm$, we impose the natural transmission conditions, i.e., the continuity of the solution and of the normal flux, namely
		\begin{equation}\label{micro_natural_trans_cond}
		\begin{aligned} 
		&&u_\epsilon^\pm &=u_\epsilon ^M &&\text{on $(0,T) \times S_{*,\epsilon}^\pm$},\\
		&&D^\pm \nabla u_\epsilon^\pm\cdot\nu &= \epsilon^\gamma D^M \left(\frac{x}{\epsilon}\right)\nabla u_\epsilon^M\cdot\nu &&\text{on $(0,T) \times S_{*,\epsilon}^\pm$}.
		\end{aligned}
		\end{equation}
	\end{subequations}
	We emphasize that we consider scalar diffusion coefficients which are constant in the bulk-domains just for an easier notation. The results can be easily extended to more general problems, e.g., to matrix-valued diffusion coefficients having an oscillating structure also in the bulk regions.
	\begin{remark}
		In order to keep the notation as clear as possible, we skip the parameter $\gamma \in [-1, 1)$ in the labeling of the solution $u_\epsilon$. The same holds for limit of $u_\epsilon$ for $\epsilon \to 0$.
	\end{remark}

	\subsubsection{Assumptions on the data:}
	\begin{itemize}
		\item[(A1)]	
		The diffusion coefficient $D_\epsilon: \Oe \rightarrow \R$ is given by
		\begin{flalign}
		D_\epsilon(x)=
		\begin{cases}
		D^\pm, &\text{ $x\in \Omega_\epsilon^\pm$,} \\
		D^M\left(\frac{x}{\epsilon}\right), &\text{ $x\in \OseM$,}\\
		\end{cases}
		\end{flalign}
		such that $D^\pm >0$ and $D^M \in L^\infty(Z^*)$, periodically extended with respect to $\bar y$ with period $Y$, and there exists $c_0>0$ such that $D^M(y) \geq c_0$, for a.e. $y \in Z^*$. 
		\item[(A2)] For the reaction rates in the bulk-domains, we assume $ f \in L^2((0,T)\times \Omega)$.
		\item[(A3)] 
		We assume $g_\epsilon \in  L^2((0,T) \times \OseM)$ with $\frac{1}{\sqrt{\epsilon}} \lVert g_\epsilon \rVert _{L^2((0,T) \times \OseM)} \leq C$, and there exists $g_0 \in  L^2((0,T) \times \Sigma \times Z^*)$ such that
		$$
		g_\epsilon \to g_0 \quad \quad \mbox{in the two-scale sense}.
		$$
		For the definition of the two-scale convergence in the channels, see Section \ref{Two_scale_convergence}.
		\item[(A4)] 
		We assume $h_\epsilon \in  L^2((0,T) \times N_\epsilon)$ with $ \lVert h_\epsilon \rVert _{L^2((0,T) \times N_\epsilon)} \leq C$, and there exists $h_0 \in  L^2((0,T) \times \Sigma \times N)$ such that 
		$$
		h_\epsilon \to h_0 \quad \quad \mbox{in the two-scale sense on }  N_\epsilon.
		$$
		Again, see Section \ref{Two_scale_convergence} for the definition of two-scale convergence on the surface $N_\epsilon$ of the channels.
		\item[(A5)] For the initial condition, we assume that 
		\begin{equation*} 
		u_{i,\epsilon} (x)=  \left\{
		\begin{array}{ll}
		u_i^\pm(x),& \text {  $x\in \Omega_\epsilon^\pm$,}\\
		u_i^M(\bar x, \frac{x}{\epsilon}),& \text {  $x \in \OseM$,}\\
		\end{array} 
		\right.
		\end{equation*}
		with $(u_i^+, u_i^M, u_i^-) \in L^2(\Omega^+)\times L^2(\Sigma, C(\bar {Z^*}))\times L^2(\Omega^-)$, and $u_i^M$ is periodic with respect to $\bar y$ with period $Y$. Especially, Lemma \ref{oscillation_lemma} implies
		$$
		\frac{1}{\sqrt{\epsilon}} \left \lVert u^M_i \left(\cdot_{\bar x}, \frac{\cdot_{x}}{\epsilon}\right) \right \rVert _{L^2 (\OseM)} \leq C ||u_i^M||_{L^2(\Sigma, C(\bar{Z^*}))} \leq C,
		$$
		and 
		$$
		u^M_i \left(\cdot_{\bar x}, \frac{\cdot_{x}}{\epsilon}\right) \to u_i^M(\bar x, y) \quad \quad \mbox{in the two-scale sense}.
		$$
		The last statement follows from \eqref{oscillation_lemma_2} in Lemma 
		\ref{oscillation_lemma}, applied to $u^M_i (\cdot_{\bar x}, \frac{\cdot_x}{\epsilon}) \psi$, where $\psi$ is a test-function in the definition of the two-scale convergence for channels.
	\end{itemize}
	\begin{example}
		\begin{itemize}
			\item[(i)] An example of a reaction rate $g_\epsilon$ satisfying assumption (A3) is given by $g_\epsilon(t,x) = g(t, \bar x, \frac{x}{\epsilon})$ where $g \in L^2((0,T) \times \Sigma, C( \overline{Z^*}))$ periodically extended with respect to $\bar y$ with period $Y$. This follows again by Lemma \ref{oscillation_lemma}.
			\item[(ii)] An example of a reaction rate $h_\epsilon$ satisfying assumption (A4) is given by $h_\epsilon(t,x) = h(t, \bar x, \frac{x}{\epsilon})$ where $h \in L^2((0,T), C( \overline{\Sigma} \times \overline{N}))$ periodically extended  with respect to $\bar y $ with period $Y$. This follows by Lemma \ref{oscillation_lemma_N}.
		\end{itemize}
	\end{example}
	
	\subsubsection{Variational formulation of the microscopic problem}
	Due to the scaling of the problem in the channel-domain, we consider function spaces with scaled inner products as follows: Let $L_\epsilon$ denote the space 
	$$
	L_\epsilon:= L^2 (\Omega_\epsilon)=L^2(\Omega_\epsilon^+) \times L^2(\OseM) \times L^2(\Omega_\epsilon^-)
	$$ 
	equipped with the inner product 
	\begin{equation*}
	(u_\epsilon, v_\epsilon)_{L_\epsilon} := \int_{\Omega_\epsilon^+} u_\epsilon v_\epsilon \,dx +\int_{\Omega_\epsilon^-} u_\epsilon v_\epsilon \,dx + \frac{1}{\epsilon} \int _{\OseM} u_\epsilon v_\epsilon  \,dx,
	\end{equation*}
	and for $\gamma \in [-1,1)$, let $H_{\gamma,\epsilon}$ be the space 
	$$
	H_{\gamma,\epsilon}: = H^1 (\Omega_\epsilon) = \left\{ (u_\epsilon^+, u_\epsilon^M, u_\epsilon^-) \in H^1(\Omega_\epsilon^+) \times H^1(\OseM) \times H^1(\Omega_\epsilon^-): u_\epsilon^\pm = u_\epsilon^M \mbox{ on } S^\pm_{*, \epsilon} \right\}
	$$
	equipped with the inner product
	\begin{equation*}
	(u_\epsilon,v_\epsilon) _ {H_{\gamma,\epsilon}} := (u_\epsilon,v_\epsilon)_{L_\epsilon} + \int _{\Omega_\epsilon^+}\nabla u_\epsilon  \nabla v_\epsilon  \,dx +\int _{\Omega_\epsilon^-}\nabla u_\epsilon  \nabla v_\epsilon  \,dx  + \epsilon^\gamma \int_{\OseM} \nabla u_\epsilon \nabla v_\epsilon \,dx.
	\end{equation*}
	The associated norms are denoted with $||\cdot||_{L_\epsilon}$ and $||\cdot||_{H_{\gamma,\epsilon}}$.
	Let $H^\prime_{\gamma,\epsilon}$ denote the dual of $H_{\gamma,\epsilon}$. We obtain the Gelfand-triple
	\begin{equation} \label{Gelfand}
	H_{\gamma,\epsilon} \hookrightarrow L_\epsilon \hookrightarrow H^\prime_{\gamma,\epsilon}.
	\end{equation}
	
	The variational formulation of the problem \eqref{microscopic_model} is given as follows. Find 
	$$
	u_\epsilon \in L^{2}((0,T), H_{\gamma,\epsilon}) \cap H^1((0,T), H^\prime_{\gamma,\epsilon})
	$$ 
	such that $u_\epsilon$ satisfies
	\begin{eqnarray} \label{weak_formulation_microsc_pb}
	&& \left<\partial _{t} u_\epsilon,\phi\right>_{H^\prime_{\gamma,\epsilon} ,H_{\gamma,\epsilon}}  + \sum_{\pm}\int_{\Omega_\epsilon^\pm} D^\pm \nabla u^\pm_\epsilon \nabla \phi\,dx  + \,  \epsilon^\gamma  \int_{\OseM}  D^M\left(\frac{x}{\epsilon}\right)\nabla u^M_\epsilon \nabla \phi \,dx \nonumber \\
	&& = \sum_{\pm}\int_{\Omega_\epsilon^\pm} f^\pm(t,x)\phi \,dx + \, \frac{1}{\epsilon}\int_{\OseM} g_\epsilon(t, x)\phi \,dx  -\int_{N_\epsilon} h_\epsilon(t,x) \phi \,dS(x),
	\end{eqnarray}
	for all $\phi \in H_{\gamma,\epsilon}$ and  a.e. $ t \in (0,T)$, together with the initial condition \eqref{micro_initial_cond_pb_1}.
	
	The duality pairing in \eqref{weak_formulation_microsc_pb} is connected to the scaling in the microscopic equation \eqref{micro_reac_diff_eq} in the following way: For $\partial_t u_\epsilon^\pm \in L^2((0,T), H^1(\Omega_\epsilon^\pm)')$ and
	$\partial_t u_\epsilon^M \in L^2((0,T), H^1(\OseM)')$, we obtain from the Gelfand-tripel \eqref{Gelfand} and the definition of the inner product in $L_\epsilon$ 
	$$
	\left<\partial _{t} u_\epsilon,\phi\right>_{H^\prime_{\gamma,\epsilon}, H_{\gamma,\epsilon}} = \sum_{\pm} \left<\partial _{t} u_\epsilon^\pm,\phi\right>_{H^1(\Omega_\epsilon^\pm)', H^1(\Omega_\epsilon^\pm)} + \frac{1}{\epsilon}\left<\partial _{t} u_\epsilon^M,\phi\right>_{H^1(\OseM)', H^1(\OseM)}.
	$$ 
	However, in our case, the time derivative $\partial_t u_\epsilon$ is a functional defined on the whole space $H_{\gamma,\epsilon}$, i.e. a space of functions defined on the whole domain $\Omega_\epsilon$, and it is not straightforward to restrict such functionals to $H^1(\Omega_\epsilon^\pm)'$ and $H^1(\OseM)'$.
	
	\section{A priori estimates}
	\label{A_priori_estimates}
	
	\begin{lemma}\label{trace_est}
		For all $\theta >0$ there exists $C(\theta) >0$ such that for all  $v_\epsilon \in  H^1(\OseM)$, we have 
		\begin{equation}\label{trace_estimate}
		\left	\lVert v_\epsilon \right\rVert_{L^2(N_\epsilon)} \leq  \frac{C(\theta)}{\sqrt{\epsilon}} \left\lVert v_\epsilon \right\rVert _{L^2(\OseM)}+ \theta \sqrt{\epsilon}\left \lVert \nabla v_\epsilon\right \rVert _{L^2(\OseM)}.
		\end{equation}
	\end{lemma}
	
	\noindent This result is easily obtained by decomposing $\OseM$ into microscopic cells $\epsilon Z_{\bar k}^*$ for $\bar k\in I_\epsilon$ together with a scaling argument and the usual trace estimate on $\partial Z^*$.
	\begin{theorem} \label{thm_apriori_est}
		The weak solution $u_\epsilon$ of the microscopic problem \eqref{weak_formulation_microsc_pb} satisfies the following a priori estimate with a constant $C >0$ independent of $\epsilon$
		\begin{eqnarray*} 
			\max_{0\leq t \leq T} \lVert u_{\epsilon}(t) \rVert _{L_\epsilon} 
			+\lVert  u_{\epsilon} \rVert _{L^2((0,T), H_{\gamma, \epsilon})} 
			+ \lVert \partial_{t} u_\epsilon \rVert _{L^2((0,T), H^\prime_{\gamma, \epsilon})}\leq C.
		\end{eqnarray*}
	\end{theorem}
	\begin{proof}
		The proof is rather standard, however, we have to exhibit the explicit dependence on the scale parameter $\epsilon$. Inserting $u_\epsilon$ as a test function in  \eqref{weak_formulation_microsc_pb}, for a.e. $t \in (0,T)$, we obtain 
		\begin{eqnarray*}
			&&\left<\partial _{t} u_\epsilon,u_\epsilon\right>_{H^\prime_{\gamma,\epsilon} ,H_{\gamma,\epsilon}}  + \sum_{\pm}\int _{\Omega_\epsilon^\pm} D^\pm \nabla u^\pm_{\epsilon} \nabla u^\pm_{\epsilon}\,dx + \ \epsilon ^ \gamma \int _{\OseM} D^M\left(\frac{x}{\epsilon}\right)\nabla u^M_{\epsilon}  \nabla u^M_{\epsilon}\,dx\\
			&& = \sum_{\pm}\int _{\Omega_\epsilon^\pm} f^\pm u^\pm_{\epsilon}\,dx+ \frac{1}{\epsilon} \int _{\OseM}g_\epsilon u^M_{\epsilon}\,dx -\int_{N_\epsilon} h _\epsilon u^M_\epsilon \,dS (x) \, \leq \, \sum_{\pm} \lVert u^\pm_{\epsilon}  \rVert^2 _{L^2 (\Omega_\epsilon^\pm) } \\
			&&+ \frac{1}{\epsilon}\left \lVert u^M_{\epsilon} \right \rVert^2 _{L^2 (\OseM) } + \left \lVert u^M_{\epsilon} \right \rVert^2 _{L^2 (N_\epsilon) } + \sum_{\pm} \lVert f^\pm \rVert^2 _{L^2 (\Omega_\epsilon^\pm) }+  \frac{1}{\epsilon}\lVert g_\epsilon \rVert^2 _{L^2 (\OseM) }+\lVert h_\epsilon\rVert^2 _{L^2(N_\epsilon)}.
		\end{eqnarray*}
		Now, the assumption (A1) on the diffusion coefficients and the trace estimate \eqref{trace_estimate} yield 
		\begin{eqnarray*}\label{a_priori_est_1}
			&& \frac{d}{dt} \lVert u_\epsilon (t) \rVert ^2_{L_\epsilon}+ \sum_{\pm} {\lVert \nabla u^\pm_{\epsilon} \rVert^2 _{L^2(\Omega_\epsilon^\pm)}} +  \epsilon ^\gamma\left\lVert \nabla u_\epsilon^M \right\rVert^2 _{L^2(\OseM)} \leq (1 +C(\theta)) \lVert u_\epsilon (t) \rVert ^2_{L_\epsilon}\\
			&&+ \theta \epsilon \left \lVert \nabla u_\epsilon^M \right \rVert ^2 _{L^2(\OseM)} + \sum_{\pm} \lVert f^\pm \rVert^2 _{L^2 (\Omega_\epsilon^\pm) }+  \frac{1}{\epsilon}\lVert g_\epsilon \rVert^2 _{L^2 (\OseM) }+\lVert h_\epsilon\rVert^2 _{L^2(N_\epsilon)}.
		\end{eqnarray*}
		For $\theta >0$ sufficiently small, the gradient term on the right hand side can be absorbed on the left-hand side. Integrating with respect to time, using assumptions (A2)-(A4) and Gronwall's inequality, we obtain 
		$$
		\max_{0\leq t \leq T} \lVert u_{\epsilon}(t) \rVert _{L_\epsilon} 
		+\lVert  u_{\epsilon} \rVert _{L^2(0,T;H_{\gamma, \epsilon})} \leq C.
		$$
		The estimate for the norm $\lVert \partial _{t} u_\epsilon \rVert _{L^2(0,T; H^\prime_{\gamma , \epsilon})} $ follows by testing \eqref{weak_formulation_microsc_pb} with  $v \in H_{\gamma, \epsilon}$ with $\lVert v \rVert _{H _{\gamma, \epsilon} }\leq 1$ and using similar arguments as above. 
	\end{proof}						
	
	\section{Two-scale convergence for channels}
	\label{Two_scale_convergence}
	
	\begin{definition} 
		\begin{itemize}
			\item[(a)]We say the sequence $v_\epsilon \in L^2((0,T) \times \OseM)$ converges in the two-scale sense to a limit function $v_0(t,\bar{x},y) \in  L^2((0,T) \times \Sigma \times Z^*)$, if we have
			\begin{equation*}
			\lim_{\epsilon \rightarrow 0} \frac{1}{\epsilon} \int_{0}^{T} \int _{\OseM} v_\epsilon(t,x) \psi \left(t,\bar{x},\frac{x}{\epsilon}\right) \,dx \,dt= \int_{0}^{T} \int _{\Sigma} \int _{Z^*} v_0(t,\bar{x},y)\psi(t,\bar{x},y) \,dy \,d\bar{x} \,dt,
			\end{equation*} 
			for all $\psi \in  L^2((0,T) \times \Sigma, C( \overline{Z^*}))$ and extended periodically with respect to $\bar{y}$ with period $Y$.
			\item[(b)]We say the sequence $v_\epsilon \in L^2((0,T) \times N_\epsilon)$ converges in the two-scale sense on $N_\epsilon$ to a limit function $v_0(t,\bar{x},y) \in  L^2((0,T) \times \Sigma \times N)$, if we have
			\begin{equation*}
			\lim_{\epsilon \rightarrow 0} \int_{0}^{T} \int _{N_\epsilon} v_\epsilon(t,x) \psi \left(t,\bar{x},\frac{x}{\epsilon}\right) \,dS(x) \,dt= \int_{0}^{T} \int _{\Sigma} \int _{N} v_0(t,\bar{x},y)\psi(t,\bar{x},y) \,dS(y) \,d\bar{x} \,dt,
			\end{equation*} 
			for all $\psi \in  L^2((0,T), C(\overline{ \Sigma} \times \overline{N}))$ and extended periodically with respect to $\bar{y}$ with period $Y$.
		\end{itemize}
	\end{definition}
	
	The following Theorem \ref{two-scale_compactness_result} provides compactness results for two-scale convergence for channels and on channels' boundaries. To prove it, we use the following two lemmas.  A similar result for bulk-domains can be found e.g. in \cite[Theorem 2]{Wall02}. Our proof is based on the decomposition of $\Sigma$ into $\epsilon$-cells which is also used in the definition of the unfolding operator, see also \cite[Lemma 4.3]{Neu07}. 
	\begin{lemma}[Oscillation-lemma for channels]\label{oscillation_lemma}
		Let $\psi \in L^1((0,T) \times \Sigma, C( \overline{Z^*}))$ and extended periodically with respect to $\bar{y}$ with period $Y$. Then $\psi\left(\cdot, \cdot,\frac{\cdot}{\epsilon}\right) \in L^1((0,T) \times \OseM)$ with 
		\begin{equation}\label{oscillation_lemma_1}
		\frac{1}{\epsilon |Z|}\left|\left|\psi\left(\cdot, \cdot,\frac{\cdot}{\epsilon}\right) \right|\right|_{L^1((0,T) \times \OseM)} \leq ||\psi||_{L^1((0,T) \times \Sigma,  C( \overline{Z^*}))}
		\end{equation}
		and
		\begin{equation}\label{oscillation_lemma_2}
		\lim_{\epsilon \rightarrow 0} \frac{1}{\epsilon} \int_{0}^{T} \int _{\OseM} \psi \left(t,\bar{x},\frac{x}{\epsilon}\right) \,dx \,dt= \int_{0}^{T} \int _{\Sigma} \int _{Z^*} \psi(t,\bar{x},y) \,dy \,d\bar{x} \,dt.
		\end{equation} 
	\end{lemma}
	\begin{proof}
		Estimate \eqref{oscillation_lemma_1} is straightforward. To prove \eqref{oscillation_lemma_2}, let us first consider $\psi \in C([0,T] \times \overline{\Sigma} \times \overline{Z^*}))$. 
		By similar arguments like in \cite[Lemma 4.3]{Neu07}, we obtain
		\begin{eqnarray*}
			\frac{1}{\epsilon} \int_{0}^{T} \int _{\OseM}  \psi \left(t,\bar{x},\frac{x}{\epsilon}\right) \,dx \,dt= \int_{0}^{T} \int_{\Sigma} \int_{Z^*} \psi \left(t, \epsilon\left[\frac{\bar x}{\epsilon} \right]+ \epsilon \bar y, y\right) \,dy\,d \bar x \,dt,
		\end{eqnarray*}
		where $[\cdot]$ denotes the Gau{\ss} bracket.
		Taking now the limit $\epsilon \to 0$, the assertion \eqref{oscillation_lemma_2} follows by Lebesgue's theorem. 
		For a general $\psi \in L^1((0,T) \times \Sigma, C( \overline{Z^*}))$, 
		the estimate \eqref{oscillation_lemma_2} follows by approximation due to \eqref{oscillation_lemma_1}.
	\end{proof}
	
	\begin{lemma}[Oscillation-lemma for channels' lateral boundaries]\label{oscillation_lemma_N}
		Let $\psi \in L^1((0,T), C(\overline{\Sigma} \times \overline{N}))$ and extended periodically with respect to $\bar{y}$ with period $Y$. Then $\psi\left(\cdot, \cdot,\frac{\cdot}{\epsilon}\right) \in L^1((0,T) \times N_\epsilon)$ with 
		\begin{equation}\label{oscillation_lemma_N_1}
		\frac{1}{|\Sigma| |N|}\left|\left|\psi\left(\cdot, \cdot,\frac{\cdot}{\epsilon}\right) \right|\right|_{L^1((0,T) \times N_\epsilon)} \leq ||\psi||_{L^1((0,T), C( \overline{ \Sigma} \times \overline{N}))}
		\end{equation}
		and
		\begin{equation}\label{oscillation_lemma_N_2}
		\lim_{\epsilon \rightarrow 0} \int_{0}^{T} \int _{N_\epsilon} \psi \left(t,\bar{x},\frac{x}{\epsilon}\right) \,dS(x) \,dt= \int_{0}^{T} \int _{\Sigma} \int _{N} \psi(t,\bar{x},y) \,dS(y) \,d\bar{x} \,dt.
		\end{equation} 
	\end{lemma}
	\begin{proof}
		The proof follows the lines of Lemma \ref{oscillation_lemma} and is skipped here. We emphasize that for the measure of $N_\epsilon$ it holds $|N_\epsilon| = O(1)$, thus no scaling of integrals on $N_\epsilon$ is needed. 
	\end{proof}

	\begin{theorem}\label{two-scale_compactness_result}
		\begin{itemize}
			\item[(i)] Let $v_\epsilon$ be a sequence of functions in $L^2((0,T) \times \OseM)$ such that 
			\begin{equation} \label{assump_comptc_result}
			\frac{1}{\sqrt{\epsilon}} \lVert v_\epsilon \rVert _{L^2((0,T) \times \OseM)} \leq C. 
			\end{equation}
			Then there exists a function $u_0(t,\bar{x},y) \in  L^2((0,T) \times \Sigma \times Z^*)$ such that, up to a sub-sequence, $v_\epsilon$ converges to $v_0$ in the  two-scale sense.
			\item[(ii)] Let $v_\epsilon$ be a sequence of functions in $L^2((0,T) \times N_\epsilon)$ such that 
			\begin{equation} \label{assump_comptc_result_N}
			\lVert v_\epsilon \rVert _{L^2((0,T) \times N_\epsilon)} \leq C. 
			\end{equation}
			Then there exists a function $v_0(t,\bar{x},y) \in  L^2((0,T) \times \Sigma \times N)$ such that, up to a sub-sequence, $v_\epsilon$ converges to $v_0$ in the  two-scale sense on $N_\epsilon$.
		\end{itemize}
	\end{theorem}
	\begin{proof}
    Using Lemma \ref{oscillation_lemma} and Lemma \ref{oscillation_lemma_N}, this result can be proved analogously to \cite[Proposition 4.2]{Neu07}, see also \cite[Theorem 1.2]{All92}.
    	\end{proof}
	
	\section{Convergence results}
	\label{Convergence_results}
	In this section, based on the a priori estimates, we prove convergence results for the weak solutions of the microscopic problems in the limit $\epsilon \to 0$. These results are the basis for the derivation of the homogenized model in the next section. 
	\begin{proposition}\label{prop_convergence_bulk}
		Let $u_\epsilon$ be the sequence of solutions of problem \eqref{weak_formulation_microsc_pb}. Then, there exists $u_0^\pm \in L^2((0,T), H^1(\Omega^{\pm}))$ such that up to a subsequence
		\begin{eqnarray*}
			& \chi_{\Omega_\epsilon^{\pm}}u_\epsilon^\pm \rightarrow u_0^\pm & \text{ strongly in } L^2((0,T) \times \Omega^{\pm}),\\
			&u_\epsilon^\pm (\cdot_t, \cdot_{\bar x}, \pm \epsilon) \rightarrow u_0^\pm|_{\Sigma} & \text{ strongly in } L^2((0,T) \times \Sigma),\\
			&\chi_{\Omega_\epsilon^{\pm}}\nabla u_\epsilon^{\pm} \rightharpoonup \nabla u_0^\pm & \text{ weakly in } L^2((0,T) \times \Omega^{\pm}).
		\end{eqnarray*}
	\end{proposition}
	\begin{proof}
		A similar result was proved in \cite[Proposition 2.1 and 2.2]{Neu07} for time derivatives $\partial_t u_\epsilon^\pm$ in  $L^2((0,T) \times \Omega_\epsilon^\pm)$. In our case, we have that $\partial_t u_\epsilon^\pm$ are functionals on the spaces $H^1_0(\Omega_\epsilon^\pm)$. However, the methods from \cite{Neu07} can easily be extended to our setting and we skip the details. 
	\end{proof}
	
	\begin{theorem} \label{thm_convergence_channels}
		Let $u_\epsilon$ be the sequence of solutions of problem \eqref{weak_formulation_microsc_pb}. Then, there exist $u_0^M \in L^2((0,T)\times \Sigma)$ and $u_1^M \in L^2 ((0,T) \times \Sigma ; H^1(Z^*) / \R )$ such that up to a subsequence
		\begin{eqnarray*}
			&u_\epsilon^M \rightarrow u_0^M & \text{ in the two-scale sense},\\
			&\epsilon^{\frac{\gamma + 1}{2}}\nabla u_\epsilon ^M \rightarrow \nabla_y u_1^M & \text{ in the two-scale sense.}
		\end{eqnarray*}
		Furthermore, for $\gamma = -1$ it holds:
		\begin{eqnarray*}
			&\nabla u_\epsilon^M \rightarrow 0 & \text{ in the two-scale sense}.
		\end{eqnarray*}
	\end{theorem}
	
	\begin{proof}
		From the a priori estimates in Theorem \ref{thm_apriori_est} we have 
		$$
		\frac{1}{\sqrt{\epsilon}}||u_\epsilon^M||_{L^2((0,T)\times\OseM)} + \epsilon^{\frac{\gamma}{2}} || \nabla u_\epsilon^M||_{L^2((0,T) \times \OseM)}  \leq C.
		$$
		Thus, by Theorem \ref{two-scale_compactness_result}, there exist functions $u_0^M \in L^2((0,T)\times \Sigma \times Z^*)$ and $\xi \in L^2((0,T)\times \Sigma \times Z^*)^n$ such that up to a subsequence 
		\begin{eqnarray*}
			&u_\epsilon^M \rightarrow u_0^M & \text{ in the two-scale sense},\\
			&\epsilon^{\frac{\gamma + 1}{2}}\nabla u_\epsilon ^M \rightarrow \xi & \text{ in the two-scale sense.}
		\end{eqnarray*}
		Next, we show that $u_0^M$ does not depend on $y \in Z^*$. Let $\psi (t, \bar{x},y) \in C^\infty_0 ((0,T) \times \Sigma \times Z^*)^n$ and periodically extended with respect to $\bar{y}$ with period $Y$. Since $\gamma <1$, the two-scale convergence of $\epsilon^{\frac{\gamma + 1}{2}}\nabla u_\epsilon ^M$ yields
		\begin{eqnarray}
		&&\lim_{\epsilon \rightarrow 0} \int_{0}^{T} \int_{\OseM} \nabla u_\epsilon^M (t,x) \cdot\psi \left(t, \bar{x},\frac{x}{\epsilon}\right) \,dx\,dt 
		=0. \label{lim_grad_ueM}
		\end{eqnarray}
		Integration by parts in \eqref{lim_grad_ueM} gives
		\begin{eqnarray*}
			0
			&=&-\lim_{\epsilon \rightarrow 0} \int_{0}^{T} \int_{\OseM} u_\epsilon^M (t,x) \left(\nabla _{\bar{x}} \cdot  \psi \left(t, \bar{x},\frac{x}{\epsilon}\right)+ \frac{1}{\epsilon}\nabla_y \cdot \psi \left(t, \bar{x},\frac{x}{\epsilon}\right)\right)\,dx\,dt \\ 	 	
			&=& - \int_{0}^{T} \int_{\Sigma} \int_{Z^*} u_0^M(t, \bar{x},y) \nabla_y \cdot \psi (t, \bar{x},y) \,dy \,d \bar{x} \,dt \label{indep_of_Z^*_1}.
		\end{eqnarray*}
		Since $Z^*$ is connected, we conclude that $u_0^M$ does not depend on $y \in Z^*$.
		
		To show that the limit $\xi$ can be represented by the gradient of a function $u_1 \in L^2 ((0,T) \times \Sigma ; H^1(Z^*) / \R)$, we consider a test-function $\psi$ as above, with the additional property that $\nabla_y \cdot  \psi =0$. Now integration by parts yields
		\begin{eqnarray}
		&&\frac{1}{\epsilon} \int_{0}^{T} \int_{\OseM} \epsilon^{\frac{\gamma + 1}{2}} \nabla u_\epsilon^M(t,x) \cdot \psi \left(t, \bar{x}, \frac{x}{\epsilon}\right) \,dx \,dt 
		\nonumber \\ && 
		= - \frac{1}{\epsilon} \int_{0}^{T} \int _{\OseM} \epsilon^{\frac{\gamma + 1}{2}} u _\epsilon^M (t,x)  \nabla  _{\bar{x}} \cdot \psi \left(t, \bar{x}, \frac{x}{\epsilon}\right) \,dx \,dt. \label{repres_limit_u1_1}
		\end{eqnarray}
		By passing to the limit $\epsilon \rightarrow 0$ in \eqref{repres_limit_u1_1}, for $\gamma \in (-1,1)$ we immediately obtain 
		\begin{equation} \label{repres_grad_1} 
		\int_{0}^{T} \int _{\Sigma} \int_{Z^*} \xi (t, \bar{x},y) \cdot \psi (t, \bar{x},y) \,dy \,d\bar{x} \,dt = 0. 
		\end{equation}
		For $\gamma =-1$ the factor $\epsilon^{\frac{\gamma + 1}{2}} =1$ in the right-hand side of \eqref{repres_limit_u1_1} gives no contribution, thus a more refined argument is necessary. Passing to the limit in  \eqref{repres_limit_u1_1} leads to 
		\begin{equation} \label{repres_grad_2} 
		\int_{0}^{T} \int _{\Sigma} \int_{Z^*} \xi (t, \bar{x},y) \cdot \psi (t, \bar{x},y) \,dy \,d\bar{x} \,dt 
		= -\int_{0}^{T} \int _{\Sigma} u_0 (t,\bar{x})\nabla_{\bar x} \cdot \int _{Z^*} \psi (t, \bar{x}, y) \,dy \, d \bar{x} \,dt.
		\end{equation}
		For the components $\psi_i $ of the test-function $\psi$ we have for $1 \leq i \leq n$:
		\begin{equation*}
		\int_{Z^*} \psi_i (t, \bar{x},y) \,dy = \int _{Z^*} \psi \cdot e_i \,dy = \int _{Z^*} \psi \cdot \nabla_{y} y_i \,dy = - \int_{Z^*} y_i\nabla_y \cdot \psi  \,dy =0. 
		\end{equation*}
		It follows that \eqref{repres_grad_1} is satisfied also for $\gamma=-1$. From \eqref{repres_grad_1}, by using similar arguments like e.g. in \cite[Chapter IV, Lemma 1.4.1]{Soh01}), we conclude that for all $\gamma \in [-1,1)$, there exists a function $u_1^M \in L^2 ((0,T) \times \Sigma ; H^1(Z^*) / \R )$ such that 
		\begin{equation*}\label{Grad_conv_zero_1}
		\xi (t, \bar{x},y)= \nabla _{y} u_1^M (t, \bar{x},y).
		\end{equation*}
		
		Now, let us show the second part of the theorem, namely that for $\gamma =-1$, $\nabla u_\epsilon^M$ converges to zero in the two-scale sense, i.e., that $ \nabla _{y} u_1^M (t, \bar{x},y)$ is, in fact, $0$. The idea is to test the  weak formulation \eqref{weak_formulation_microsc_pb} by a suitable test-function and to pass to the limit for $\epsilon \to 0$. Let $ \alpha \in C^\infty_0 (0,T), \beta \in C^\infty _0 (\Sigma)$ and $\phi \in C^\infty (\overline{Z^*})$. We extend $\phi$ periodically with respect to $\bar y$ with period $Y$. Furthermore, we extend $\phi |_{S_*^{\pm}}$ to a function $\tilde \phi^\pm \in W^{1, \infty}(S^\pm)$ and then periodically with respect to $\bar y$ with period $Y$ such that we obtain $\tilde{\phi}^\pm(\cdot, \pm 1) \in W^{1, \infty}(\R^{n-1})$. Since $S_*^{\pm}$ is Lipschitz, such an extension exists, see e.g., \cite[Theorem 9.7]{Bre11}. 
		Finally, let $\rho \in C^\infty ([0, \infty) )$ such that $ 0 \leq \rho \leq 1, \rho = 1$ in $[0,1]$ and $\rho = 0$ in $[2, \infty)$. Now, we define $ \phi_\epsilon$ on $(0,T) \times \Omega_\epsilon$ by 
		\begin{flalign*}
		\phi_\epsilon (t,x)=
		\begin{cases}
		\epsilon \alpha (t) \beta (\bar{x}) \tilde {\phi}^+ \left (\frac{\bar{x}}{\epsilon},1\right) \rho \left(\frac{x_n}{\epsilon}\right), & (t,x)\in (0,T) \times \Omega_\epsilon^+, \\
		\epsilon  \alpha (t) \beta (\bar{x})\phi \left(\frac{x}{\epsilon}\right),& (t,x) \in (0,T) \times \OseM, \\
		\epsilon  \alpha (t) \beta (\bar{x})\tilde{\phi}^-  \left(\frac{\bar{x}}{\epsilon}, -1\right)\rho \left(-\frac{x_n}{\epsilon}\right),&  (t,x)\in (0,T) \times \Omega_\epsilon^-.
		\end{cases}
		\end{flalign*}
		By construction, $\phi _\epsilon \in L^2 (0,T; H_{\gamma, \epsilon} )$ and using it as a test-function in \eqref{weak_formulation_microsc_pb} yields: 
		\begin{eqnarray} \label{weak_formulation_microsc_pb_phi_epsilon}
		&& \left<\partial _{t} u_\epsilon,\phi_\epsilon\right>_{H^\prime_{\gamma,\epsilon} ,H_{\gamma,\epsilon}}  
		+ \sum_\pm \epsilon \int_{\Omega_\epsilon^\pm} D^\pm \alpha (t)\nabla u^\pm_\epsilon \nabla \beta (\bar{x}) \tilde {\phi}^\pm \left (\frac{\bar{x}}{\epsilon},\pm1\right) \rho \left(\pm\frac{x_n}{\epsilon}\right) \,dx \nonumber \\
		&& + \sum_{\pm}  \int_{\Omega_\epsilon^\pm} D^\pm \alpha (t) \nabla u^\pm_\epsilon  \beta (\bar{x}) \nabla_y \tilde {\phi}^\pm \left (\frac{\bar{x}}{\epsilon},\pm1\right) \rho \left(\pm \frac{x_n}{\epsilon}\right) \,dx \nonumber\\ 
		&& + \sum_{\pm}  \int_{\Omega_\epsilon^\pm} D^\pm \alpha (t) \nabla u^\pm_\epsilon  \beta (\bar{x})  \tilde {\phi}^\pm \left (\frac{\bar{x}}{\epsilon},\pm1\right) \nabla_y \rho \left(\pm \frac{x_n}{\epsilon}\right) \,dx \nonumber\\ 
		&& +  \int_{\OseM}   D^M\left(\frac{x}{\epsilon}\right) \alpha (t) \nabla u^M_\epsilon \nabla \left( \beta (\bar{x}) \phi \left (\frac{x}{\epsilon}\right) \right) \,dx  \nonumber\\
		&&=  \sum_{\pm} \epsilon \int_{\Omega_\epsilon^\pm} f^\pm(t,x)\alpha (t) \beta (\bar{x}) \tilde {\phi}^\pm \left (\frac{\bar{x}}{\epsilon},\pm1\right) \rho \left(\pm\frac{x_n}{\epsilon}\right) \,dx \nonumber \\
		&&+ \ \int_{\OseM} g_\epsilon(t, x) \alpha (t) \beta (\bar{x}) \phi \left (\frac{x}{\epsilon}\right) \,dx   - \epsilon \int_{N_\epsilon}  h_\epsilon(t, x)  \alpha (t) \beta (\bar{x}) \phi \left (\frac{x}{\epsilon}\right) \,dS(x). \,
		\end{eqnarray}
		It turns out that after integration with respect to time in \eqref{weak_formulation_microsc_pb_phi_epsilon}, in the limit $\epsilon \to 0$ all terms tend to zero, apart form the term
		$$
		\int_0^T \int_{\OseM}   D^M\left(\frac{x}{\epsilon}\right) \alpha (t) \nabla u^M_\epsilon \nabla \left( \beta (\bar{x}) \phi \left (\frac{x}{\epsilon}\right) \right) \,dx\, dt.
		$$
		This follows from the convergence properties of the sequence $u_\epsilon$ from Proposition \ref{prop_convergence_bulk} and the first part of this theorem, from the assumptions (A2)-(A4) for the reaction rates and from the properties of the test-function $\phi_\epsilon$.  Let us illustrate this in the following. For the term involving the time derivative, we have
		\begin{equation*}
		\int_{0}^{T}<\partial_t u _\epsilon , \phi_\epsilon > _{H_{\gamma, \epsilon}^\prime , H_{\gamma, \epsilon}} \,dt	= 
		- \int_{0}^{T} (\partial _{t} \phi_\epsilon, u_\epsilon )_{L_\epsilon} \,dt 
		\end {equation*}
		and the right-hand side converges to zero due to the convergence properties of the sequence $u_\epsilon$. Terms involving $\rho$ can be estimated as follows:
		\begin{eqnarray*}
		&& \left|\int_{0}^{T}\int _{\Omega_\epsilon^+} D^+ \alpha (t) \frac{\partial u_\epsilon^+}{\partial x _n} \beta (\bar{x}) \tilde{\phi}^+ \left ( \frac{\bar{x}}{\epsilon},1\right )  \rho ^ \prime \left(\frac{x_n}{\epsilon}\right)  \,dx \,dt \right| \\
		&& \leq C \left \lVert \frac{\partial u^+ _{\epsilon}}{\partial x _n} \right \rVert _{L^2(0,T; L^2 (\Omega_\epsilon ^+))} \left[\int_{\epsilon} ^{2\epsilon} \left | \rho^ \prime \left (\frac{x_n}{\epsilon}\right) \right |^2  \,dx_n \right]^{\frac{1}{2}}  \leq  C  \sqrt{\epsilon}\left[\int_{1} ^{2} \left | \rho^ \prime  (\lambda) \right |^2  \,d \lambda  \right]^{\frac{1}{2}}.
		\end{eqnarray*}
		Consequently, also this term converges to zero. Thus, taking the limit $\epsilon \to 0$ in \eqref{weak_formulation_microsc_pb_phi_epsilon} yields 
		\begin{equation}\label{repres_limit_u_1}
		\int_{0}^{T} \int_{\Sigma} \int_{Z^*} D^M (y)\alpha (t)  \beta (\bar{x})  \nabla _{y} \phi (y)  \nabla _{y } u_1^M (t, \bar{x},y) \,dy \,d\bar{x} \,dt =0,
		\end{equation}
		for all $(\alpha, \beta, \phi) \in C^\infty_0 (0,T) \times C^\infty _0 (\Sigma) \times C^\infty (\overline{Z^*})$. Thus, we have 
		\begin{equation}\label{2-repres_limit_u_1}
		\int_{Z^*}  D^M(y) \nabla _{y} \phi (y)  \nabla _{y } u_1^M (t, \bar{x},y) \,dy =0, \, \text{for a.e. $(t, \bar{x}) \in (0,T) \times \Sigma$}.
		\end{equation}
		By a density argument, we can insert $\phi = u_1^M(t, \bar x,.)$ in \eqref{2-repres_limit_u_1}, which concludes the proof.
		\end{proof}
		
		In the last theorem of this section, we establish relations across the interface $\Sigma$ between the limit $u_0^M$ in the channels and the limits $u_0^\pm$ in the bulk-domains. These relations contribute to the transmission conditions in the homogenized model.
		\begin{theorem}
		Let $\gamma \in [-1, 1)$. The following continuity relations  across $\Sigma$ hold between the limit functions $u_0^\pm$ from Proposition \ref{prop_convergence_bulk} and $u_0^M$ from Theorem \ref{thm_convergence_channels}: 
		\begin{equation}\label{continuity_macroscopic_sol}
		u_0^M (t, \bar x) = u_0^\pm (t, \bar x,0) \quad \mbox{ for a.e. } (t, \bar x) \in (0,T) \times \Sigma,
		\end{equation}
		especially, we have $u_0^+ |_{\Sigma}(t, \bar x)  = u_0^- |_\Sigma (t, \bar x)$ a.e. on $(0,T) \times \Sigma $.
		\end{theorem}
		\begin{proof}
		Let $\psi (t, \bar{x},y) \in C^\infty ((0,T) \times \Sigma \times Z^*)^n$ such that $\psi(t,\bar x, \cdot)$ has compact support in $Z^* \cup S_*^+ \cup S_*^-$ and extend this function by zero to $Z$ and then Y-periodically with respect to $\bar{y}$. As in \eqref{lim_grad_ueM}, we have 
		\begin{equation}
		\lim_{\epsilon \rightarrow 0} \int_{0}^{T} \int_{\OseM} \nabla u_\epsilon^M (t,x) \cdot\psi \left(t, \bar{x},\frac{x}{\epsilon}\right) \,dx\,dt  = 0.
		\end{equation}
		Integrating by parts and using the continuity of the microscopic solution on $S_{*,\epsilon}^\pm$, the two-scale convergence of $u_\epsilon^M$ and the strong convergence of $u_\epsilon^\pm (\cdot_t, \cdot_{\bar x}, \pm \epsilon)$ which implies the two-scale convergence of $u_\epsilon^\pm(\cdot_t, \cdot_{\bar x}, \pm \epsilon)$ on $S^\pm_\epsilon$, we obtain:
		\begin{eqnarray*}
		0 &=& \lim_{\epsilon \to 0} - \frac{1}{\epsilon} \int_{0}^{T} \int_{\OseM} u_\epsilon^M (t,x) \left( \epsilon \nabla_{\bar x} \cdot \psi \left(t, \bar{x},\frac{x}{\epsilon}\right) + \nabla_y \cdot \psi \left(t, \bar{x},\frac{x}{\epsilon}\right) \right)\,dx\,dt \\
		&- &  \lim_{\epsilon \to 0} \, \sum_{\pm} \int_{0}^{T} \int_{S_{*,\epsilon}^\pm} u_\epsilon^\pm (t,x)  \psi \left(t, \bar{x},\frac{x}{\epsilon}\right)\cdot \nu\,dx\,dt \\
		&=&-  \int_0^T \int_\Sigma \int_{Z^*} u_0^M(t, \bar x) \nabla_y \cdot \psi(t, \bar x, y) dy\, d\bar x \, dt \\
		&- & \sum_{\pm} \int_{0}^{T} \int_\Sigma \int_{S_{*}^\pm} u_0^\pm (t,\bar x, 0)  \psi \left(t, \bar{x},y\right)\cdot \nu \,dy \,d \bar x\,dt.
		\end{eqnarray*}
		Using again the integration by parts and taking into account that $u_0^M$ does not depend on $y$ gives the desired result. 
		\end{proof}
		
		\section{Derivation of the macroscopic model}
		\label{Derivation_macroscopic_model}
		
		For the study of the macroscopic problem, we consider the space $L^2_*$ given by
		$$
		L^2_*:= L^2(\Omega^+) \times L^2(\Sigma) \times L^2(\Omega^-)
		$$
		together with the inner product
		$$
		(u_0, v_0)_{L^2_*} = \int_{\Omega^+} u_0^+ v_0^+ dx + |Z^*|\int_\Sigma u_0|_\Sigma \, v_0|_\Sigma d\bar x+ \int_{\Omega^-} u_0^- v_0^- dx,
		$$
		and we use the following Gelfand-triple
		\begin{equation}\label{Gelfand_macro}
		H^1(\Omega)  \hookrightarrow L^2_* \hookrightarrow H^1(\Omega)'.
		\end{equation}
		
		\begin{theorem}\label{macroscopic_model}
		Let $\gamma\in [-1, 1)$. The limit function $u_0= (u_0^+, u_0^M, u_0^-)$ from Proposition \ref{prop_convergence_bulk} and Theorem \ref{thm_convergence_channels} satisfies
		$u_0 \in  L^2((0,T), H^1(\Omega)) \cap H^1((0,T), H^1(\Omega)')$ with respect to the Gelfand-triple \eqref{Gelfand_macro}, 
		and is the unique weak solution of the following problem:
		\begin{equation*}\label{macro_reac_diff_eq}
		\begin{aligned} 
		&&\partial_{t} u_0^\pm -D^\pm \Delta u_0^\pm  &= f^\pm && \mbox{ in }  (0,T)\times  \Omega^\pm,  \\
		&&-D^\pm \nabla u_0^\pm \cdot \nu^\pm &= 0 && \mbox{ on }  (0,T)\times (\partial \Omega^\pm \setminus \Sigma),\\
		&& u_0^- &= u_0^+= u_0^M && \mbox{ on }  (0,T)\times \Sigma, \\
		&& [D \nabla u_0 \cdot \nu]_\Sigma = & - |Z^*|\partial_{t}u_0^M  + \int_{Z^*}g_0(\cdot_t, \cdot_{\bar x}, y) dy - \int_N h_0(\cdot_t, \cdot_{\bar x}, y) dS(y) &&\text{ on $(0,T)\times \Sigma$},\\
		&& u_0^\pm (0) & = u_i^\pm && \mbox{ in } \Omega^\pm,\\
		&& u_0^M(0) &=\frac{1}{|Z^*|} \int_{Z^*} u_i^M (\cdot_{\bar x}, y) dy&& \mbox{ on } \Sigma,
		\end{aligned}
		\end{equation*}
		where $[D \nabla u_0 \cdot \nu]_\Sigma := (D^+ \nabla u_0^+ - D^- \nabla u_0^-) \cdot \nu^+$ denotes the jump in the normal flux of the homogenized solution $u_0$ at the interface $\Sigma$.
	\end{theorem}
	\begin{proof}
		Let $\psi \in C_0^\infty([0,T) \times \overline \Omega)$ and consider as a test-function in \eqref{weak_formulation_microsc_pb} 
		$$
		\psi_\epsilon(t, x) = \left\{ \begin{array}{ll} \psi(t, \bar x, \frac{H}{H-\epsilon}(x_n \mp \epsilon)), & x \in \Omega_\epsilon^\pm,\\
		\psi(t, \bar x, 0), & x \in \OseM.
		\end{array} \right.
		$$
		Integration with respect to time, integration by parts in time, Proposition \ref{prop_convergence_bulk} and Theorem \ref{thm_convergence_channels} imply for $\epsilon \to 0$ the variational equality
		\begin{eqnarray}\label{eq_macroscopic_var}
		&&\sum_{\pm} \left\{ - \int_{(0,T) \times \Omega^\pm} u^\pm_0 \partial_t \psi^\pm  dx dt + \int_{(0,T) \times \Omega^\pm} D^\pm \nabla u_0^\pm \nabla \psi^\pm \, dx dt \right\} \nonumber\\
		&&- \, |Z^*| \int_{(0,T) \times \Sigma} u_0^M \partial _t \psi |_\Sigma \, d\bar x dt \nonumber\\
		&&= \sum_{\pm} \left\{ \int_{(0,T) \times \Omega^\pm} f^\pm \psi^\pm dx dt + \int_{\Omega^\pm} u_i^\pm  \psi^\pm(0) dx \right\} \\
		&&+  \int_{(0,T) \times \Sigma} \int_{Z^*} g_0 (t,\bar x, y) dy \, \psi |_{\Sigma}(t, \bar x) \, d\bar x dt \nonumber\\
		&& -   \int_{(0,T) \times \Sigma} \int_{N} h_0 (t,\bar x, y) dS(y) \, \psi |_{\Sigma}(t, \bar x) d\bar x dt +  \int_{\Sigma} \int_{Z^*} u_i^M ( \bar x, y) dy \, \psi |_{\Sigma}(0, \bar x) d\bar x, \nonumber
		\end{eqnarray}
		for all $\psi \in C_0^\infty([0,T) \times \overline \Omega)$. From \eqref{eq_macroscopic_var} it follows that $\partial_t u_0 \in L^2((0,T), H^1(\Omega)')$ with respect to the Gelfand-triple \eqref{Gelfand_macro} and thus, \eqref{eq_macroscopic_var} together with the continuity relation \eqref{continuity_macroscopic_sol} is precisely the weak formulation of the macroscopic problem from the Theorem. Especially, it follows that
		$
		u_0(0)= \left(u_i^+, \frac{1}{|Z^*|} \int_{Z^*} u_i^M ( \cdot_{\bar x}, y) dy, u_i^-\right).
		$
		The uniqueness of the weak solution follows by standard arguments.
	\end{proof}
	
	\begin{remark}
		We remark that the homogenized model does not depend on the parameter $\gamma \in [-1, 1)$. This is different from the result in \cite{Gah17}, where a thin layer with a heterogeneous structure is considered, with coefficients uniformly elliptic within the whole layer. There, different results were obtained for $\gamma = -1$ and $\gamma \in (-1,1)$. The result in the present paper is similar to the result obtained in \cite{Gah17} for $\gamma \in (-1,1)$. The case $\gamma = -1$ turned out to be a critical case in \cite{Gah17}, where a reaction-diffusion equation on $\Sigma$ in form of a dynamic Wentzell interface condition is obtained. The appearance of such an interface condition in the case of a thin layer perforated by channels is hindered by the disconnected channels which don't allow a macroscopic diffusion within the interface $\Sigma$.
	\end{remark}
	
	\section*{Acknowledgements}
	AB acknowledges the support by the RTG 2339 IntComSin of the German Research Foundation (DFG). MG was supported by the SCIDATOS project, funded by the Klaus Tschira Foundation (Grant number 00.0277.2015), and the Odysseus program of the Research Foundation - Flanders FWO (Project-Nr. G0G1316N). Parts of the paper have been created during a research stay of AB at the Hasselt University.

\end{document}